\documentclass[11pt]{article}
\usepackage[english]{babel}
\usepackage[latin1]{inputenc}
\usepackage{exscale}
\usepackage[centertags]{amsmath}
\usepackage{amsfonts,amstext,amssymb,amsthm}
\usepackage{comment}
\usepackage{newlfont}
\usepackage{enumerate}
\usepackage{graphicx}
\usepackage{these}
\usepackage{color}
\usepackage{hyperref}
\usepackage{url}
\title{The Isoperimetric Profile of a Noncompact Riemannian Manifold for Small Volumes}
\author{Stefano Nardulli}
\begin{document}
      \maketitle
      \tableofcontents       
      \addcontentsline{toc}{section}{\numberline{}Introduction}
      \newpage
\section{Introduction}
Let $(M, g)$ be an $n$-dimensional Riemannian manifold. We deal mainly with the problem of finding a relatively compact domain $D\subset\subset M$ that minimizes $Area(\partial D)$ among domains of the same volume, for sufficiently small values of volume. We reformulate the problem in the language of currents from geometric measure theory. 
Given $0<v<Vol(M)$, consider all integral currents $T$ in $M$ with volume $v$, and denote the mass of the boundary as $Area(\partial T)$. From now on we consider the problem of finding minimizing currents with a fixed volume constraint. This problem is referred as the \textbf{isoperimetric problem} throughout the paper. 

When we speak about area and volume, respectively $Area(\cdot)$ and $Vol(\cdot)$, we do not mention the metric when this is clear from the context, but sometimes it will be necessary to specify the metric for the sake of clarity and according to this convention we may write $Area_g$ and $Vol_g$ where $g$ denotes the metric. 
 
The principal achievements of this paper concern the link between the theory of pseudo-bubbles and the isoperimetric problem for small volumes in a complete Riemannian manifold with some kind of "boundedness at infinity" on the metric and its fourth derivatives. This task was carried out by the same author in the context of manifolds for which there is existence of minimizers in all volumes, in particular for manifolds with co-compact isometry group or manifolds with finite volume (compare with \cite{RRosales}). In this paper, we deal with the same questions, but with entirely new techniques, which meet the difficulties arising from the lack of existence of minimizers. Namely, we isometrically embed the manifold $M$ into a metric space which is composed of a disjoint union of pieces $(M_{\infty},p_{\infty},g_{\infty})$.  These pieces are limits manifolds of sequences $(M,p_j,g)_j$, with $p_j\in M$, in some suitable pointed $C^{k,\alpha}$ topology. The arguments presented here are useful because they allow us to prove non-trivial phenomena for $M$ complete, non-compact, possibly without existence of minimizers, provided that sufficiently many sequences $(M,p_j,g)$ have a limit in a $C^{k,\alpha}$ topology. For the convenience of the reader, we repeat the relevant material from \cite{NarAnn},\cite{BM}, \cite{Pet}, \cite{PXu} and \cite{NarGD} without proofs for a self-contained exposition.\\ 

First we recall the definition of a pseudo-bubble. Let $Q=id-P$, where $P$ is orthogonal projection of $L^2 (T^1 _p M)$ on the first eigenspace of the Laplacian $T^1 _p M$ is the fiber over $p$ of the unit tangent bundle of the Riemannian manifold $M$. 
\begin{Def}\cite{NarAnn}
          A \textbf{pseudo-bubble} is an hypersurface $\mathcal{N}$ embedded in $M$ such that there exists a point $p\in M$ 
          and a function $u$ belonging to $C^{2,\alpha}(T^1_p M\backsimeq\mathbb{S}^{n-1},\mathbb{R})$, such that $\mathcal{N}$ is the graph of
$u$ in normal polar coordinates centered at $p$, i.e. $\mathcal{N}=\left\lbrace exp_p (u(\theta )\theta ),\:\theta\in T^1_p M \right\rbrace$ and $Q(H(u))$ is a real constant, where $H$ is the mean curvature operator. 
\end{Def}\noindent
To state a uniqueness theorem for pseudo-bubbles we need the notion of \textit{center of mass}.
\begin{Def}
Let $(\Omega ,\mu )$ be a probability space and $f:\Omega\rightarrow M$ a measurable function. We consider the following function $\mathcal{E}:M\rightarrow [0,+\infty[$:
$$\mathcal{E}(x):=\frac{1}{2}\int_{\Omega } d^2 (x,f(y))d\mu (y).$$
The \textbf{center of mass} of $f$ with respect to the measure $\mu$ is the minimum of 
$\mathcal{E}$ on $M$, provided that it exists and is unique.
\end{Def}\noindent
In particular, we can speak about the center of mass of a hypersurface of small diameter (we apply this definition to the $(n-1)$-dimensional measure of the boundary).
The main result on pseudo-bubbles is the following theorem.
\begin{Thm}[\cite{NarAnn}, Theorem $1$]\label{TT1}
Let $M$  be a complete Riemannian manifold. Let $\mathcal{F}^{k,\alpha}$ be the fiber bundle on $M$ whose fiber over $p$ is the space of $C^{k,\alpha}$ functions on the unit tangent sphere $T^1_p M$. There exists a $C^{\infty}$ map, $\beta:M\times ]0,Vol(M)[\rightarrow\mathcal{F}^{2,\alpha}$ such that for all $p\in M$, and all sufficiently small $v>0$, the hypersurface $exp_p (\beta(p,v)(\theta )\theta)$ is the unique pseudo-bubble with center of mass $p$ enclosing a volume $v$.
\end{Thm}\noindent
\paragraph{Remark:} In general the volume $v$ depends on the point $p$, but it can be chosen to depend continuously on the point $p$. In the case of manifolds with $C^{2,\alpha}$-bounded geometry, there is a uniform upper bound $\bar{v}=\bar{v}(n,k,\alpha, Q)$ such that the conclusion of the preceding theorem is true for every $p\in M$ and every $0<v<\bar{v}$.
\paragraph{Remark:} If $g$ is an isometry of $M$, $g$ sends pseudo-bubbles to pseudo-bubbles and $g\circ\beta =\beta\circ g$ ($g$ acts only on the first factor $M$).\\

\subsection{Main Results}
According to \cite{MJ}, small solutions of the isoperimetric problem in compact Riemannian manifolds, or noncompact manifolds with cocompact isometry group, are close to geodesic balls. Namely, they are graphs in normal coordinates of $C^{2,\alpha}$ small functions. This holds also for non-compact manifolds under a $C^4$ bounded geometry assumption, which will be proven in section $3$. In any case, it follows that these small isoperimetric domains are pseudo-bubbles. 
\paragraph{Remark:} $C^4$ boundedness is due only to the technical limits of the methods employed for proving theorem \ref{Regthm}. 
The main result of this paper is theorem \ref{1}, which provides a criterion for existence of minimizers having sufficiently small volume. To state this theorem correctly, let us recall the basic definitions from the theory of convergence of manifolds, as articulated in \cite{Pet}. 
\begin{Def} [Petersen \cite{Pet}] A sequence of pointed complete Riemannian manifolds is said to converge in
the pointed $C^{m,\alpha}$ topology $(M_i, p_i, g_i)\rightarrow (M,p,g)$ if for every $R > 0$ we can find
a domain $\Omega_R$ with $B(p,R)\subseteq\Omega_R\subseteq M$, a natural number $\nu_R\in\mathbb{N}$, and embeddings $F_{i,R}:\Omega_R\rightarrow M_i$ for $i\geq\nu_R$ such that
$B(p_i,R)\subseteq F_{i,R} (\Omega_R)$ and $F_{i,R}^*(g_i)\rightarrow g$ on $\Omega_R$ in the $C^{m,\alpha}$ topology.
\end{Def}\noindent 
It is easy to see that this type of convergence implies pointed Gromov-Hausdorff convergence, because $C^{0,\alpha}$ convergence implies Lipschitz convergence and Lipschitz convergence implies Gromov-Hausdorff convergence, see \cite{msrnrs} theorem 3.7 page 74. When all manifolds in question are closed, then the maps $F_i$ are diffeomorphisms. So for closed manifolds we can speak about unpointed convergence. In this case, convergence can only happen if all the manifolds in the tail end of the sequence are diffeomorphic. In particular, classes of closed Riemannian manifolds that are precompact in some $C^{m,\alpha}$ topology contain at most finitely many diffeomorphism types. For the precise definition of $C^{m,\alpha}$ bounded geometry, see the definition below.
\begin{Def} [Petersen \cite{Pet}]
Suppose $A$ is a subset of a Riemannian $n$-manifold
$(M,g)$. We say that the $C^{m,\alpha}$-norm on the scale of $r$ of $A\subseteq(M,g)$ satisfies
$||A||_{C^{m,\alpha},r}\leq Q$ if we can find charts
$\psi_s :\mathbb{R}^n\supseteq B(0,r)\rightarrow U_s\subseteq M$ such that
\begin{enumerate}[(i):]
      \item For all $p\in A$ there exists $U_s$ such that $B(p,\frac{1}{10}e^{-Q}r)\subseteq U_s$.
      \item $|D\psi_s|\leq e^Q$ on $B(0,r)$ and $|D\psi_s^{-1}|\leq e^Q$ on $U_s$.
      \item $r^{|j|+\alpha}||D^jg_s||_{\alpha}\leq Q$
 for all multi indices j with $0\leq |j|\leq m$, where $g_s$ is the matrix of functions of metric coefficients in the $\psi_s$ coordinates regarded as a matrix on $B(0,r)$. 
\end{enumerate}      
\end{Def}
\begin{Def} For fixed $Q>0$, $n\geq 2$, $m\geq 0$, $\alpha\in ]0,1]$, and $r>0$, define $\mathcal{M}^{m,\alpha}(n, Q, r)$ as the class of complete, 
pointed Riemannian $n$-manifolds $(M,p,g)$ with $||M||_{C^{m,\alpha},r}\leq Q$. 
\end{Def}
In the sequel, $n\geq 2$, $r,Q>0$, $m\geq 4$, $\alpha\in [0,1]$. 
\begin{Res}\label{1}
There exists $0<v^*=v^*(n,r,Q,m,\alpha)$ such that for all $M\in\mathcal{M}^{m,\alpha}(n, Q, r)$ and for every $v$ such that $0<v<v^*$ then 
\begin{enumerate}[(I):]          
           \item \label{Main0-} The two following statements are equivalent,
           \begin{enumerate}[(a):] 
           \item \label{MainI} the function $p\mapsto f_M(p,v)$ attains its minimum, 
           \item \label{MainII} there exists solutions of the isoperimetric problem at volume $v$,
           \end{enumerate}
            \item \label{Main0} $I_M(v)=Min\{f_{M_{\infty}}(p_{\infty},v)|\:(M,p_j, g)\rightarrow (M_{\infty},p_{\infty},g)\: for\: some\: (p_j)\}$.
\end{enumerate}
Here $p_j\in M$ and the function $p\mapsto f_M(p,v)$ gives the area of pseudo-bubbles contained in a given manifold $M$, with center of mass $p\in M$ and enclosed volume $v$. Moreover, every solution $D$ of the isoperimetric problem is s.t. $\partial D=\left\lbrace exp_{p_0} (u(\theta )\theta ),\:\theta\in T^1_{p_0} M \right\rbrace$, where $p_0$ is a minimum of $p\mapsto f_M(p,v)$ and conversely. With $\beta$ obtained in theorem \ref{TT1}. $f_M$ is invariant and $\beta$ equivariant under the group of isometries of $M$.
\end{Res}\noindent
The proof of theorem \ref{1} will be achieved at the end of section \ref{s3}.\\ 

\textbf{Remark:} The interest in theorem \ref{1} is the reduction of the infinite-dimensional problem of finding a minimizer to a finite dimensional one, namely, to the problem of finding the minima of a smooth function defined on the manifold $M$.\\

Let us mention one important consequence (theorem \ref{Cor1}) of the isoperimetric profile defined below. 
\begin{Def}
     Let $M$ be a Riemannian manifold of dimension $n$ (possibly with infinite volume).
     Denote by $\tau_{M}$ the set of  relatively compact open subsets of $M$ with smooth boundary. The function $I:[0,Vol(M)[\rightarrow [0,+\infty [$ such that
      $I(0)=0$
     \Fonct{I}{]0,Vol(M)[}{[0, +\infty [}{v}{Inf_{\left\lbrace \begin{array}{l}
                                                                                                       \Omega\in \tau_{M}\\
                                                                                                       Vol(\Omega )=v 
                                                                                               \end{array} \right\rbrace }  \{ Area(\partial \Omega )\} } 
is called the \textbf{isoperimetric profile function} (or shortly the \textbf{isoperimetric profile}) of the manifold $M$.
\end{Def}\noindent
In this respect, we need to compute an asymptotic expansion of the function $v\mapsto f(p,v)$. For this we use results of \cite{PXu}. For the sake of completeness, the statement of the following theorem is included. Note that any term denoted $\mathcal{O} ({r^k})$ here is a smooth function on $\mathbb{S}^{n-1}$ that might depend on $p$ but which is bounded by a constant independent of $p$ times $r^k$ in the $C^2 $ topology.
\begin{Def}
We denote by $c_n := \frac{Area(\mathbb{S}^{n-1})}{[Vol (\mathbb{B}^{n})]^{\frac{n-1}{n}}}$ the constant in the Euclidean isoperimetric profile. 
\end{Def}  
The following lemma uses calculation partially computed in \cite{PXu}.    
\begin{Lemme}[\cite{NarAnn}]\label{PXvol} \label{AeApb} Asymptotic expansion of the area of pseudo-bubbles as a function of the enclosed volume.  
 \begin{equation}\label{Ipexp}
f(p,v)=c_n v^{\frac{n-1}{n}}\left\lbrace 1+a_p \left(\frac{v}{\omega_n} \right) ^{\frac{2}{n}} +\mathcal{O} (v^{\frac{4}{n}}) \right\rbrace,  
\end{equation} 
with $a_p:=-\frac{1}{2n(n+2)}Sc(p)$.                        
\end{Lemme}
Denote by $Sc$ the scalar curvature function of $M$.
\begin{Res} \label{Cor1} For all $M\in\mathcal{M}^{m,\alpha}(n, Q, r)$, let $$S=Sup_{p\in M}\{ Sc(p)\}.$$
                Then the isoperimetric profile $I_{M}(v) $ has the following asymptotic expansion in a neighborhood of the origin
                             \begin{equation}\label{fasymex}
I_{M}(v)=c_n v^{\frac{n-1}{n}}\left( 1-\frac{S}{2n(n+2)}\left( \frac{v}{\omega_n}\right) ^{\frac{2}{n}} +o(v^{\frac{2}{n}})\right),
\end{equation}
where $o(t^{\alpha})$ indicates a function $g:]-\varepsilon, \varepsilon[\rightarrow\mathbb{R}$, with $\varepsilon >0$, such that $lim_{t\rightarrow 0}\frac{g(t)}{t^{\alpha}}=0$.
\end{Res}
In theorem \ref{Cor1} and lemma \ref{AeApb}, $\mathcal{O}(t^{\alpha})$ and $o(t^{\alpha})$ are functions that depend only on $t$.
The asymptotic expansion of the volume of pseudo-bubbles and the volume of their boundary can be computed with lemma \ref{PXvol}, which yields an expansion for the profile. 
    
\subsection{Plan of the article}
\begin{enumerate}
           \item Section $2$ describes why and in what sense approximate solutions of the isoperimetric problem, in the case of small volumes, are close to Euclidean balls, providing a decomposition theorem for domains belonging to an almost minimizing sequences in small volumes. 
           \item In section $3$ we prove theorem \ref{1}, generalizing to the case of $C^4$-bounded geometry manifolds some results of \cite{NarAnn}, in particular corollary \ref{Kleiner1} which constitutes the only known proof to my knowledge of the fact that for small volumes minimizers are invariant under the action of the groups of isometries of $M$ that fix their barycenters.
           \item In section $4$ the results of preceding sections and those of \cite{NarGD}, \cite{MJ}, \cite{PXu} are applied to obtain the first two non-zero coefficients in the asymptotic expansion of the isoperimetric profile in the non-compact case under $C^4$-bounded geometry assumption on $M$.
           \end{enumerate}               
\subsection{Acknowledgements}  
I would to thank my former Ph.D. advisor Pierre Pansu for attracting my attention to the subject of this paper, and for his comments, suggestions and encouragement which helped shape this article. I am also indebted to Frank Morgan, Manuel Ritor\'e, and Efstratios Vernadakis for their useful comments and remarks. Finally, I would like to express my gratitude to Renata Grimaldi for numerous mathematical discussions. 
\newpage
      \section{Partitions of domains}
\subsection{Introduction}
In this section it is assumed that   
\begin{enumerate}
           \item $M$ has bounded geometry ($|\mathcal{K}|\leq\Lambda$ and $inj_{M}\geq\varepsilon >0$) where $inj_{M}$ is the injectivity radius of $M$ ,
           \item the domains $D_j \in\tau_{M}$ are \textbf{approximate solutions} i.e. $\frac{Area(\partial D_j)}{I(Vol (D_j ))}\rightarrow 1$ for $j\rightarrow +\infty$.
\end{enumerate} 

We prove in this section the following theorem.

\begin{Res}\label{tprec1}
Let $(M,g)$ be a Riemannian manifold with bounded geometry, $D_j $ a sequence of approximate solutions of the isoperimetric problem such that $Vol_g (D_j )\rightarrow 0$. Then there exist $p_j\in M$, and radii $R_j\rightarrow 0$ such that  
\begin{equation}
            \lim_{j\rightarrow +\infty}\frac{Vol(D_j \Delta B(p_j , R_j ))}{Vol(D_j )}\rightarrow 0.
\end{equation}
\end{Res}     
The proof of theorem \ref{tprec1} occupies the rest of this section.
\subsection{Euclidean version of theorem \ref{tprec1}}
Roughly speaking, we have that in $\mathbb{R}^n$ approximate solutions of the isoperimetric problem are close to balls in the mass norm, as stated in the following theorem.
A good reference for  this result is  \cite{LeoRigot}.
\begin{Thm}\label{JT4}
           Let $\left\{ T_j \right\}\subset \mathbb{I}_n (\mathbb{R}^n )$ be a sequence of integral currents, satisfying
     $$\lim_{j\rightarrow +\infty}\frac{\mathbf{M}(\partial T_j )}{\mathbf{M}(T_j )^{\frac{n-1}{n}}}=c_n .$$
           Then there exist balls $W_j$ such that up to a subsequence 
                    $$\frac{\mathbf{M}(T_j \Delta W_j)}{\mathbf{M}(W_j )}\rightarrow 0.$$
\end{Thm}
\begin{Sdem}
            We can use here the theory of $BV$ functions and that of finite perimeter sets as described in \cite{Giu} because for all polyhedral chains $P$, $||\chi_{Spt||P||}||_{BV(\mathbb{R}^n )}<+\infty$. In what follows we translate our problem into the language of $BV$ functions.\\ \indent
 Let $|\cdot |$ be the Lebesgue measure on $\mathbb{R}^n$. Now we give an argument regarding minimizing sequences that will be useful in the sequel. Let $(E_k)_{k\geq 1}$ be a minimizing sequence of domains for the functional $\mathcal{H}^{n-1}(\partial (\cdot ))$ such that $|E_k|=1$.\\ \\
1. A compactness theorem stated in \cite{Giu} (page 17) ensures that there exists a set $E$ such that a subsequence   
$$\chi_{E_k}\rightarrow \chi_{E}$$  
in $L^1 _{loc}(\mathbb{R}^n )$.\\
2. By lower semicontinuity of Lebesgue measure and of the perimeter function, it follows that $|E|\leq \liminf_{k\rightarrow +\infty} |E_k |\leq 1$, $$P(E,\mathbb{R}^n )\leq  \liminf_{k\rightarrow +\infty} P(E_k ,\mathbb{R}^n )\leq c_n.$$
Now if we show that $|E|=1$, then we finish the proof, because Euclidean isoperimetric domains are round balls, so $E$ is the Euclidean ball of volume $1$. This together with $L^1 (B(0,2))$ convergence ensure that the mass outside this Euclidean ball goes to zero and that the volume of the set-theoretic symmetric difference $|E\Delta E_k|$ goes to zero.\\  
A clear proof that $|E|=1$ for Carnot-Caratheodory groups is given in \cite{LeoRigot}, and for this reason we will not repeat it here. It occurs in two steps: 
\begin{itemize}
        \item to show that there exist translates of $E_k$ having an intersection with the ball of radius $1$ of mass not less than a constant $m_0 >0$ (Lemma $4.1$ of  \cite{LeoRigot}), 
        \item to argue that we cannot find a non-negligible subset of $E_k$ far away from this radius $1$ ball because  $E_k$ is almost perimeter-minimizing among all sets of measure $1$ (Lemma $4.2$,  \cite{LeoRigot}). 
\end{itemize}
To prove the theorem it is sufficient to apply the preceding argument to sets $E_j$ obtained from dilating $supp ||T_j ||$  by a factor of $\frac{1}{\mathbf{M}(T_j )^{\frac{1}{n}}}$ and setting $W_j$ equal to $\mathbf{M}(T_j )^{\frac{1}{n}}E$ . 
\end{Sdem} 
\subsection{Lebesgue numbers}
Let $(M,g)$ be a Riemannian manifold with bounded geometry. We can construct a good covering of  $M$ by balls having the same radius.
\begin{Lemme}\label{Lebesgue}
Let $(M,g)$ be a Riemannian manifold with bounded geometry. There exist an integer $N$, some constants $C$, $\epsilon>0$ and a covering $\mathcal{U}$ of $M$ by balls having the same radius $3\epsilon$ and having also the following properties.
\begin{enumerate}
  \item $\epsilon$ is a Lebesgue number for $\mathcal{U}$, i.e. every ball of radius $\epsilon$ is entirely contained in at least one element of $\mathcal{U}$ and meets at most $N$ elements of $\mathcal{U}$.
  \item For every ball $B$ of this covering, there exists a $C$ bi-Lipschitz diffeomorphism on a Euclidean ball of the same radius.
\end{enumerate}
\end{Lemme}

\begin{Dem} Let $\epsilon=\frac{inj_{M}}{3}$. Let $\mathcal{B}=\{ B(p,\epsilon) \}$ be a maximal family of balls of $M$ of radius $\epsilon$ that have the property that any pair of distinct members of $\mathcal{B} $ have empty intersection. Then the family $2\mathcal{B}:=\{ B(p,2\epsilon) \}$ is a covering of $M$.
Furthermore, for all $y\in M$, there exist $B(p,\varepsilon )\in\mathcal{B}$  such that $y\in B(p,2\epsilon)$ and thus $B(y,\varepsilon )\subseteq B(p,3\epsilon)$. Hence $\epsilon$ is a Lebesgue number for the covering $3\mathcal{B}$. Let $B(p,3\epsilon)$ and $B(p',3\epsilon)$ be two balls of $3\mathcal{B}$ having nonempty intersection. Then $d(p,p')<6\epsilon$, hence $B(p',\epsilon)\subseteq B(p,7\epsilon)$. The ratios $Vol(B(p,7\epsilon))/Vol(B(p,\epsilon))$ are uniformly bounded because the Ricci curvature of $M$ is bounded from below, and hence the Bishop-Gromov inequality applies. The number $N$ of disjoint balls of radius $\epsilon$, contained in $B(p,7\epsilon)$, is bounded and does not depend on $p$. Thus the number of balls of $3\mathcal{B}$ that intersect one of them is uniformly bounded by $N$. We conclude the proof by taking $\mathcal{U}:=3\mathcal{B}$.
In fact by Rauch's comparison theorem, see Cheeger-Ebin \cite{ChEbin} page 29, for every ball $B(p,\epsilon)$, the exponential map is 
$C$ bi-Lipschitz with a constant $C$ that depends only on $\epsilon$ and on upper bounds for the sectional curvature $\mathcal{K}$.
\end{Dem}

\subsection{Partition domains in small diameter subdomains} 
This section is inspired by the article of B\'erard and Meyer \cite{BM} lemma II.15 and the theorem of appendix C, page 531.
\begin{Prop}
Let $I$ be the isoperimetric profile of $M$. Then
      $$ \limsup_{a\rightarrow 0}\frac{I(a)}{a^{\frac{n-1}{n}}}\leq c_n .$$
\end{Prop}
\begin{Dem} Fix a point $p\in M$.
      \begin{eqnarray*}
            \limsup_{a\rightarrow 0}\frac{I(a)}{a^{\frac{n-1}{n}}}\leq  \limsup_{a\rightarrow 0}
            \frac{Area(\partial B(p,r(a)))}{Vol(B(p,r(a)))^{\frac{n-1}{n}}}
      \end{eqnarray*}
      with $r(a)$ such that $Vol(B(p,r(a)))=a$. Changing variables in the limits, we find
      \begin{eqnarray*}
            \limsup_{a\rightarrow 0}\frac{Area(\partial B(p,r(a)))}{Vol(B(p,r(a)))^{\frac{n-1}{n}}} & = & \limsup_{r\rightarrow 0}\frac{Area(\partial B(p,r))}{Vol(B(p,r))^{\frac{n-1}{n}}}\\
         \limsup_{r\rightarrow 0}\frac{r^{n-1}Area(\mathbb{S}^{n-1})+\cdots }{[r^{n}Vol(\mathbb{B}^{n})+\cdots ]^{\frac{n-1}{n}}}& = & c_n .      
      \end{eqnarray*}
\end{Dem}
  \begin{Def}
      Let $r>0$. We define the \emph{unit grid} of $\mathbb{R}^{n}$,  $G_1$, as the set of points which have at least one integer coordinate. We call $G$  a \emph{grid of mesh} $r$ if $G$ is of the form $x+rG_1$ where $x=(x_1,...,x_n)\in \mathbb{R}^{n}$. 
      We denote by $\mathcal{G}_r:=([0,r]^n , \mathcal{L}^n ) $ the set of all grids of mesh $r$, endowed with the natural Lebesgue measure given by the bijection $\Phi:[0,r]^n\rightarrow \mathcal{G}_r$, $\Phi:x\to x+rG_1$, on $([0,r]^n ,dx_1\cdots dx_n)$.       
\end{Def}
\begin{Prop}\label{l5}
        Let $D$ be an open subset of $\mathbb{R}^n$.
       $$ \frac{1}{r^n}\int_{\mathcal{G}_r} Area(D\cap G) \mathcal{L}^n (dG)=\frac{n}{r}Vol (D). $$
\end{Prop}
\begin{Dem}
      We observe that every grid $G$ of mesh $r$ decomposes as a union of $n$ sets $G^{(i)}$ of the type $x+rG_1^{(i)}$ where $G_1^{(i)}$ is the set of points with integer $i-$th coordinate.\\
Moreover $G^{(i)}\cap G^{(j)}$ has $(n-1)$-dimensional Hausdorff measure equal to zero, for every $i\neq j$. This latter fact allows us to ensure that
\begin{equation}
Area(D\cap G)=\sum_{i=1}^n Area(D\cap G^{(i)}).
\end{equation}
Thus 
      \begin{equation}
      \frac{1}{r^n}\int_{\mathcal{G}_r} Area(D\cap G) \mathcal{L}^n (dG) = \frac{1}{r^n}\sum_{i=1} ^n  \int_{[0,r]^n }Area(D\cap G^{(i)}) \mathcal{L}^n (dG),
            \end{equation}
            but
            \begin{equation}
            \int_{[0,r]^n }Area(D\cap G^{(i)}) \mathcal{L}^n (dG)=\int_{[0,r]^n }Area(D\cap G_x^{(i)}) dx_1\cdots dx_n
            \end{equation}
            by the identification done by $\Phi$ and 
            \begin{equation}
            \int_{[0,r]^n }Area(D\cap G_x^{(i)}) dx_1\cdots dx_n = \int_{0}^r\left\{\int_{[0,r]^{n-1}} Area(D\cap G_x^{(i)})dx_1\cdots\widehat{dx_i}\cdots dx_n\right\} dx_i
            \end{equation}
            by Fubini theorem and 
            \begin{equation}
            \int_{[0,r]^{n-1}} Area(D\cap G_x^{(i)})dx_1\cdots\widehat{dx_i}\cdots dx_n= r^{n-1} Area(D\cap G_{(0,...,x_i,...0)}^{(i)}),
\end{equation}
            by domain invariance.
            It follows 
            \begin{equation}
            \int_{[0,r]^n }Area(D\cap G_x^{(i)}) dx_1\cdots dx_n=\int_{0}^r r^{n-1} Area(D\cap G_{(0,...,x_i,...0)}^{(i)})dx_i,
            \end{equation}
            \begin{equation}
          \frac{1}{r^n}\sum_{i=1} ^n  \int_{0}^r r^{n-1} Area(D\cap G^{(i)}) \mathcal{L}^n (dG) = \frac{n}{r}Vol (D),          
      \end{equation}
       that finally gives 
       \begin{equation}
        \frac{1}{r^n}\int_{\mathcal{G}_r} Area(D\cap G) \mathcal{L}^n (dG)=\frac{n}{r}Vol (D).
       \end{equation}
\end{Dem}   
\begin{Cor}
Let $r>0$. Let $D$ be an open set of $\mathbb{R}^n$.  There exists a grid $G$ of mesh $r$ such that  
\begin{equation}
Area(D\cap G)\leq \frac{n}{r}Vol (D).
\end{equation}  
\end{Cor}   
\begin{Prop}\label{l1}
      We denote $D_{G,k}$ the connected components of $D\setminus G$. Then
      $$ \frac{\sum_k Area(\partial D_{G,k})-Area(\partial D)}{Vol(D)^{\frac{n-1}{n}}}\rightarrow 0$$
      as $\displaystyle
          \frac{Vol(D)^{\frac{1}{n}}}{r}\rightarrow 0.$ 
\end{Prop}
\begin{Dem}
      For every grid $G$,
      $$\sum_k Area(\partial D_{G,k})-Area(\partial D)=2Area(D\cap G).$$ By corollary \ref{l5}, there exists a grid $G$ such that $Area(D\cap G)\leq \frac{n}{r}Vol (D)$. We deduce that 
      \begin{eqnarray*}
             0\leq\frac{\sum_k Area(\partial D_{G,k})-Area(\partial D)}{Vol(D)^{\frac{n-1}{n}}}\leq\frac{\frac{2n}{r}Vol (D)}
             {Vol(D)^{\frac{n-1}{n}}}=\frac{2nVol (D) ^{\frac{1}{n}}}{r}.
      \end{eqnarray*}
      Thus if $r$ is very large with respect to $Vol(D)^{\frac{1}{n}}$ then  
      $$  \frac{\sum_k Area(\partial D_{G,k})-Area(\partial D)}{Vol(D)^{\frac{n-1}{n}}}$$
      is close to $0$.
\end{Dem}

\begin{Prop}\label{l2}
      Let $M$ be a Riemannian manifold with bounded geometry. Let $D_j$ be a sequence of domains of  $M$ so that 
      \begin{enumerate}
              \item $Vol (D_j )\rightarrow 0$.
    \item $\limsup_{j\rightarrow +\infty} \frac{Area(\partial D_j )}{Vol(D_j )^{\frac{n-1}{n}}}\leq c_n$.
      \end{enumerate}
      For any sequence $( r_j  )$ of positive real numbers that tends to zero ($r_j \rightarrow 0$ ) and  $\frac{Vol(D_j )^{\frac{1}{n}}}{r_j }\rightarrow 0$,
      there exists a partition $D_j =\bigcup_k D_{j,k}$ of $D_j$ in domains $D_{j,k}$ with
       $Diam(D_{j,k})\leq const_{M}\cdot r_j$ such that 
      $$\limsup_{j\rightarrow +\infty} \frac{\sum_k Area(\partial D_{j,k})}{(\sum_k Vol(D_{j,k} ))^{\frac{n-1}{n}}}\leq c_n .$$       
\end{Prop}
\begin{Dem}
We apply lemma \ref{Lebesgue} and we take a covering $\{ \mathcal{U}\}$ of $M$ by balls of radius $3\epsilon$, of multiplicity $N$, i.e., $N$ is as in lemma 
\ref{Lebesgue} and Lebesgue number $\epsilon>0$. For every ball $B(p,3\epsilon)$ of this family, we fix a diffeomorphism $\phi_p :B(p,3\epsilon)\to B_{\mathbb{R}^n}(0,3\epsilon)$ of  Lipschitz constant $C$. Observe here that by bounded geometry assumptions $C$ could be chosen in such a way it is independent of $p$. For every $j$ we fix also a radius $r_j >>Vol(D_j)^{\frac{1}{n}}$ and we map the grids of mesh $r_j$ of $\mathbb{R}^n$ in $B(p,3\epsilon)$ via $\phi_p$, i.e. for $G\in\mathcal{G}_{r_{j}}$, we have
$$G_p =\phi_{p}^{-1}(G).$$
Let us denote by $D_{j,k}$ the connected components of $D_j \setminus (\cup_p G_p )$.
We are looking for an estimate of the supplementary boundary volume introduced by the partition in this $D_{j,k}$,
$$\sum_k Area(\partial D_{j,k})-Area(\partial D_j )=2Area(D_j \cap (\cup_l G_l )).$$ 
First estimate the average $m=   \frac{1}{r_j ^n }\int_{\mathcal{G}_{r_j}}Area(D_j\cap (\cup_l G_l ))\mathcal{L}^n(dG)$ of this volume over all possible choices of the grids $G\in\mathcal{G}_{r_{j}}$.
\begin{eqnarray*}
m   & \leq & 
      \frac{1}{r_j ^n }\sum_p \int_{\mathcal{G}_{r_j}}Area(D_j\cap G_p )\mathcal{L}^n (dG) \\
      & \leq &  \frac{1}{r_j ^n }\sum_p \int_{\mathcal{G}_{r_j}}{Area}_{(\mathbb{R}^n , {\phi_p ^{-1}}^{*} (g))}(\phi_p (D_j )\cap G)\mathcal{L}^n(dG)\\
      & \leq & 
      \frac{C}{r_j ^n }\sum_p \int_{\mathcal{G}_{r_j}}{Area}_{(\mathbb{R}^n , can)}(\phi_p (D_j\cap\mathcal{U}_p  )\cap G)\mathcal{L}^n(dG)\\
      & \leq & C\frac{n}{r_j }\sum_p Vol (\phi_p (D_j \cap B(p,3\epsilon)))\\
      & \leq & C^2\frac{n}{r_j }\sum_p Vol (D_j \cap B(p,3\epsilon))\\
      & \leq & C^2 \frac{n}{r_j}NVol(D_j ).
\end{eqnarray*}
This is true because every point of $M$ is contained in at most $N$ balls $B(p,3\epsilon)$. Then there exists $G $ in $\mathcal{G}_{r_j}$ such that 
$$Area(D_j \cap (\cup_p G_p ))\leq C^2 \frac{n}{r_j}NVol(D_j ),$$ and so
$$0\leq\frac{\sum_k Area(\partial D_{j,k})-Area(\partial D_j )}{Vol(D_j )^{\frac{n-1}{n}}}\leq  2 C^2 \frac{n}{r_j}NVol(D_j )^{\frac{1}{n}}.$$
From the last inequality we obtain    
$$\limsup_{j\rightarrow +\infty} \frac{\sum_k Area^{M}(\partial D_{j,k})}{(\sum_k Vol ^{M}(D_{j,k} ))^{\frac{n-1}{n}}}\leq
\limsup_{j\rightarrow 0} \frac{Area^{M}(\partial D_j )}{Vol ^{M}(D_j )^{\frac{n-1}{n}}}\leq c_n .$$
Now, fix $x\in D_j\setminus (\cup_p G_p )$. By construction, $\epsilon$ is a Lebesgue number of the covering $\{ \mathcal{U}\}$, and there exists a ball $B(p,3\epsilon)$ that contains $B(x,\epsilon)$. Let $D_{j,k}$ denote the connected component of $D_j\setminus (\cup_p G_p )$ that contains $x$, and $D'_{j,k}$ the connected component of $\phi_p (B(p, 3\epsilon))\setminus G$ that contains $\phi_p (x)$. We observe that $D'_{j,k}$ is contained in a cube of edge $r_j$; if $j$ is large enough so that $r_j \leq \epsilon/C\sqrt{n}$, then $D'_{j,k}$ is contained in $\phi_p (B(x,\epsilon))$, hence $D_{j,k}$ is contained in $\phi_{p}^{-1}D'_{j,k}$, which has diameter at most $C\,r_j$.
\end{Dem}

\subsection{Selecting a large subdomain}

We first show that an almost Euclidean isoperimetric inequality can be applied to small domains.
\begin{Lemme}\label{lem1}Let $M$ be a Riemannian manifold with bounded geometry.
                          Then
                          \begin{equation}
                                   \frac{Area(\partial D)}{Vol(D)^{\frac{n-1}{n}}}\geq c_n (1-\eta (diam(D)))
                          \end{equation} 
with $\eta\rightarrow 0$ as $diam(D)\rightarrow 0$.
\end{Lemme}
\begin{Dem}
In a ball of radius $r<inj(M)$, we reduce to the Euclidian isoperimetric inequality via the exponential map, that is a $C$ bi-Lipschitz diffeomorphism with $C=1+\mathcal{O}(r^2 )$. This implies for all domains of diameter $<r$, 
$$ \frac{Area(\partial D)}{Vol(D)^{\frac{n-1}{n}}}\geq c_n C^{-2n+2}=c_n (1-\mathcal{O}(r^2 )).$$   
\end{Dem}

Second, we have a combinatorial lemma that tells that in a partition the largest domain contains almost all the volume.

\begin{Lemme}\label{l3}
      Let $f_{j,k}\in [0,1]$ be numbers such that for all $j $, $\sum_k f_{j,k} =1$.
      Then
      $$ \limsup_{j\rightarrow + \infty} \sum_k f_{j,k}^{\frac{n-1}{n}}\leq 1 $$
      implies that
      $$ \lim_{j\rightarrow +\infty } \max_k  f_{j,k}=1.$$
\end{Lemme}
\begin{Dem}
      We argue by contradiction.
      Suppose there exists $\varepsilon >0$ for which there exists $j_{\varepsilon}\in \mathbb{N}$ so that for all  $j\geq j_{\varepsilon}$, we have $\max_k \{f_{j,k}\}\leq 1-\varepsilon $.
      Then for all $j\geq j_{\varepsilon}$, we have $f_{j,k}\leq 1-\varepsilon $. From this inequality,
 $$\sum_k f_{j,k}^{\frac{n-1}{n}}=\sum_k f_{j,k}f_{j,k}^{\frac{-1}{n}}\geq \frac{\sum_k f_{j,k}}{(1-\varepsilon )^{\frac{1}{n}}}
        \geq \frac{1}{(1-\varepsilon)^{\frac{1}{n}}},$$ 
      hence
      $$\limsup_{j\rightarrow + \infty} \sum_k f_{j,k}^{\frac{n-1}{n}}\geq \frac{1}{(1-\varepsilon)^{\frac{1}{n}}}> 1,$$
      which is a contradiction.  
\end{Dem}
 \begin{Prop}\label{l4}
Let $M$ be a Riemannian manifold with bounded geometry. Let $D_j$ be a sequence of approximate solutions in $M$ with volumes that tend to zero. Let $r_j$ be a sequence of positive real numbers such that $r_j \rightarrow 0$ and $\frac{Vol(D_j )^{\frac{1}{n}}}{r_j }\rightarrow 0$. There exist $p_j \in M$ and $\varepsilon_j \leq const_{M}r_j$ and subdomains $D'_j\subset D_j$ such that 
\begin{enumerate}
            \item $D'_j \subseteq B(p_j ,\varepsilon_j )$
            \item $\frac{Area(\partial D'_j)}{Vol(D'_j)^{\frac{n-1}{n}}}\rightarrow c_n$
            \item $\lim_{j\rightarrow +\infty } \frac{{Vol }(D'_j )}{{Vol }(D_j )}=1.$
\end{enumerate} 
 \end{Prop}
 \begin{Dem}
       Apply proposition \ref{l2}. By the definition of isoperimetric profile and lemma \ref{lem1} we have\\
       $$Area(\partial D_{j,k})\geq I(Vol(D_{j,k} ))\geq c_n 
Vol(D_{j,k})^{\frac{n-1}{n}}(1-\eta_j )$$ where $\eta_j \rightarrow 0$. Since
  \begin{eqnarray*}
       \limsup_{j\rightarrow +\infty }\frac{\sum_k c_n 
Vol(D_{j,k})^{\frac{n-1}{n}}(1-\eta_j )}{Vol(D_j )^{\frac{n-1}{n}}} \leq
       \limsup_{j\rightarrow +\infty }\frac{\sum_k Area(\partial 
D_{j,k})}{Vol(D_j )^{\frac{n-1}{n}}} \leq c_n ,\\ 
  \end{eqnarray*}
 $$\limsup_{j\rightarrow +\infty }\frac{\sum_k 
Vol(D_{j,k})^{\frac{n-1}{n}}}{Vol(D_j )^{\frac{n-1}{n}}}\leq \limsup_{j\rightarrow +\infty }
   \frac{1}{1-\eta_j }=1.$$
 Now, set $f_{j,k}=\frac{Vol(D_{j,k})}{Vol(D_j)}$. We can suppose that $f_{j,1}=max_k \{f_{j,k}\}$. We apply lemma 
\ref{l3} and we deduce that  
 $$\frac{Vol(D_{j,1})}{Vol(D_j )}\rightarrow 1.$$
 But by construction $D_{j,1}\subset B_{M}(p_j , 
const_{M}r_j )$ for some sequence of points $p_j$ in $M$. Finally, proposition \ref{l2} gives 
$$\limsup \frac{Area(\partial D_j,1 )}{Vol(D_j )^{\frac{n-1}{n}}}\leq\limsup\frac{}{}\leq c_n .$$
Indeed, since $\{D_{j,k}\}$ is a partition of $D_j$ we have $Vol(D_j)=\sum_kVol(D_{j,k})$. Thus one can take $D'_j =D_{j,1}$.
 \end{Dem}
\subsection{End of the proof of  theorem \ref{tprec1}}
In this subsection we terminate the proof of theorem \ref{tprec1}.\\ \indent
\begin{Dem}
Let $D_j$ be a sequence of approximate solutions with $Vol(D_j )\rightarrow 0$.
According to proposition \ref{l4} there exist subdomains $D'_j \subseteq D_j$, points $p_j \in M$ and radii $\varepsilon_j \rightarrow 0$ such that 
\begin{enumerate}[(i):]
           \item $D'_j \subseteq B(p_j ,\varepsilon_j )$.
           \item $\frac{Vol(D'_j )}{Vol( D_j )}\rightarrow 1$.
           \item $\frac{Area(\partial D'_j)}{Vol(D'_j)^{\frac{n-1}{n}}}\rightarrow c_n$.
\end{enumerate} 
We identify all tangent spaces $T_{p_j }M$ with a fixed Euclidean space $\mathbb{R}^n$ and consider the domains $D''_j = exp^{-1}(D'_j )$ in $\mathbb{R}^n$. Since the pulled back metrics $\tilde{g}_j=exp_{p_j}^{*}(g_{M})$ converge to the Euclidean metric, 
$$\frac{Area(\partial D''_j)}{Vol(D''_j)^{\frac{n-1}{n}}}\rightarrow c_n. $$
According to theorem \ref{JT4}, there exist Euclidean balls $W_j =B_{eucl.}(\tilde{q}_j , R_j)$ in $\mathbb{R}^n$ such that 
$$\frac{Vol_{eucl.}(D''_j\Delta W_j)}{Vol_{eucl.}(D''_j)}\rightarrow 0.$$
Note that $\tilde {g}_j $-balls are close to Euclidean balls, thus 
$$\frac{Vol_{eucl.}(D''_j \Delta B^{\tilde{g}_j}(\tilde{q}_j ,R_j))}{Vol_{eucl.}(D''_j)}\rightarrow 0,$$
where $B^{\tilde{g}_j}(\tilde{q}_j ,R))$ is the geodesic ball of radius $R$ in $(\mathbb{R}^n, \tilde {g}_j)$, and then, for $q_j =exp_{p_j}(\tilde{q}_j)$,
$$\frac{Vol_{eucl.}(D'_j \Delta B^{g}(\tilde{q}_j ,R_j))}{Vol_{eucl.}(D'_j)}=\frac{Vol_{eucl.}(D''_j \Delta B^{\tilde{g}_j}(\tilde{q}_j ,R_j))}{Vol_{\tilde{g}}(W_j)}\rightarrow 0.$$
Finally, since $\frac{Vol(D_j \Delta D'_j)}{Vol(D_j)}\rightarrow 0$, 
$\frac{Vol_{g}(D_j\Delta B(q_j , R_j))}{Vol_{g}(D_j)}\rightarrow 0,$ \\
it is trivial to completes the proof of theorem \ref{tprec1}.  
\end{Dem}
\subsection{Case of exact solutions}

\paragraph{Remark:} When we consider the \emph{solutions} of the isoperimetric problem (this is the case treated in \cite{MJ}), and not \emph{approximate solutions}, the conclusion is stronger. In fact we can prove directly by the monotonicity formula that $D_j$ is of small diameter and this simplifies a lot the arguments showing that $D_j$ are close in flat norm to a round ball. 

\begin{Lemme}
      Assume $D_j$ is a sequence of solution of the isoperimetric problem. The dilated domains $D'''_j :=\frac{exp_{p_j}^{-1}(D_j )}{Vol_g (D_j)^{\frac{1}{n}}}$ are of bounded diameter and hence we can find a positive constant $R>0$ in the proof of the preceding theorem so that for all $j\in \mathbb{N}$ we have $$D'''_j \subseteq B(0,R).$$ 
\end{Lemme}
\begin{Dem}
      For the domains $D'''_j$, the mean curvature of the boundary in 
      $(\mathbb{R}^n,eucl)$ $h_j^{eucl}\leq M=const.$ for all $j$ (apply the L\'evy-Gromov isoperimetric inequality \cite{Gr1}, \cite{Gr2} analogously to the argument given in the proof of theorem 2.2 of \cite{MJ}) and hence the monotonicity formula of \cite{All}[5.1 (3)] page 446 gives for a fixed $r_0$ and all $j$
      \begin{equation}
      ||\partial D'''_j ||(B(a_j , r_0))\geq e^{-M r_0}\Theta^{n-1}(||\partial D'''_j||, a_j)\omega_{n-1}r_0^{n-1}
      \end{equation}
      $a_j\in spt||\partial D'''_j ||$, for a fixed $r_0$ and all $j$.
      We argue 
      $$const\geq Area_{g_{can}} (\partial D'''_j )\geq \left[ \frac{Diam_{g_{can}}(D'''_j )}{2r_0}\right]\omega_{n-1}r_0 ^{n-1}$$      
      and we can conclude that $Diam_{g_{can}}(D'''_j )$ is uniformly bounded.
\end{Dem}
\newpage
      \section{Existence for small volumes.}\label{s3}
For compact manifolds, the regularity theorem of \cite{MJ} applies, and there is no need to use the more general theorem \ref{Regthm}. For noncompact manifolds the situation is quite involved.
\subsection{Minimizers are pseudo-bubbles.}\label{321}
When $M$ is noncompact, the regularity theorem of \cite{MJ} has to be replaced by a more general statement, for the following reasons.
\begin{enumerate}
           \item Solutions of the isoperimetric problem need not exist in $M$.
           \item Minimizing sequences may escape to infinity, therefore varying ambient metrics cannot be avoided.   
\end{enumerate} 
Now, let us recall the basic result from the theory of convergence of manifolds, as exposed in \cite{Pet}. 
\begin{Thm}[Fundamental Theorem of Convergence Theory. \cite{Pet} Theorem $72$]\label{Fthmct}

$\mathcal{M}^{m,\alpha}(n, Q, r)$ is compact in the pointed $C^{m,\beta}$ topology for all $\beta <\alpha$.
\end{Thm}
In subsequent arguments will be needed a regularity theorem, in a variable metrics context.  
\begin{Thm}\label{Regthm}\cite{NarGD}
\label{T4}
Let $M^n$ be a compact Riemannian manifold, $g_j$ a sequence of Riemannian metrics of class $C^{\infty}$ that converges to a fixed metric $g_{\infty}$ in the $C^4$ topology. Assume that $B$ is a domain of  $M$ with smooth boundary $\partial B$, and $T_j$ is a sequence of currents minimizing area under volume contraints in $(M^n , g_j )$ satisfying 
$$(*): Vol_{g_{\infty}} (B\Delta T_j)\rightarrow 0.$$ 
Then $\partial T_j$ is the graph  in normal exponential coordinates of a function  $u_j$ on $\partial B$.
Furthermore, for all $\alpha\in ]0,1[$, $u_j \in C^{2,\alpha}(\partial B)$ and $||u_j||_{C^{2,\alpha}(\partial B)}\rightarrow 0$ as $j\rightarrow +\infty$.
\end{Thm}\noindent
\textbf{Remark:} Roughly speaking, theorem \ref{Regthm} says that if an integral rectifiable current $T$ is minimizing and sufficiently close in flat norm to a smooth current then $T$ is smooth too. \\
\textbf{Remark:} Theorems \ref{Fthmct} and \ref{Regthm} are the main reasons for assuming to work under $C^4$ bounded geometry assumptions in this paper.\\

In the sequel we use often the following classical isoperimetric inequality due to Pierre Berard and Daniel Meyer. 
\begin{Thm}\label{BerMeyer}(\cite{BM} Appendix C]). Let
$M^{n+1}$ be a smooth, complete Riemannian manifold, possibly with boundary, of bounded geometry (bounded sectional curvature and positive injectivity radius).
Then, given $0<\delta<1$, (the interesting case is when $\delta$ is close to $1$) there exists $v_0>0$ such that any open set $U$ of volume
$0<v<v_0$ satisfies
\begin{equation}
        Area(\partial U)\geq\delta c_n v^{\dfrac{n-1}{n}}.
\end{equation}
\end{Thm}\noindent      
\textbf{Remark:} The preceding theorem implies in particular that for a complete Riemannian manifold with bounded sectional curvature and strictly positive injectivity radius holds $I_M(v)\sim c_n v^{\frac{n-1}{n}}$ as $v\rightarrow 0$.


\begin{Lemme}\label{tpreec} Let $M\in\mathcal{M}^{m,\alpha}(n, Q, r)$, and $(D_j)$ a sequence of solutions of the isoperimetric problem with $Vol_g (D_j )\rightarrow 0$. Then possibly extracting a subsequence, there exist points $p_j\in M$ such that the domains $D_j $ are graphs in polar normal coordinates of functions $u_j$ of class $C^{2,\alpha}$ on the unit sphere of $T_{p_j} M$ of the form $u_j=r_j(1+v_j )$ with 
$||v_j||_{C^{2,\alpha }(\partial B_{T_p M}(0, 1))}\rightarrow 0$ and radii $r_j\rightarrow 0$.
\end{Lemme}
\begin{Dem}
We consider tangent spaces $T_{p_j} M$ in this situation we identify them with a fixed copy of $\mathbb{R}^n$ and in this fixed space we carry almost the same analysis as already done in \cite{NarAnn} Lemma 3.1 the only difference is that there by compactness of $M$ we can deal just with one point instead of a sequence of points $p_j$ . In fact we take domains $T_j$ to be $exp_{p_j}^{-1}(D_j)$ rescaled by $\frac{1}{r_j}$ in the same fixed copy of $\mathbb{R}^n$ then $T_j$ is a solution of the isoperimetric problem for the rescaled pulled-back metric $g_j=\frac{1}{r_j ^2}exp_p^* (g)$ which converges volumewise to a unit ball. Since the sequence $g_j$ converges at least $C^4$ to a Euclidean metric, because of the $C^4$ bounded geometry assumption on $g$ the same arguments as in the preceding lemma applies.
\end{Dem}
\begin{Lemme}\label{lm1}
For all $n,r,Q,m\geq 4,\alpha$, there exists $0<v_1=v_1(n,r,Q,m,\alpha)$ such that for all $M\in\mathcal{M}^{m,\alpha}(n, Q, r)$, for every domain $D$ solution of the isoperimetric problem with $0<Vol(D)\leq v_1$, there exists a point $p_{D}\in M$ (depending on $D$) such that $D$ is the normal graph of a function $u_D \in C^{2,\alpha}(\mathbb{S}^{n-1})$ with $u_{D}=r_{D}(1+v_{D})$ and $||v_{D}||_{C^{2,\alpha }(\mathbb{S}^{n-1})}\rightarrow 0$ as $ Vol(D)\rightarrow 0$.  
\end{Lemme}
\begin{Dem}
Otherwise there exists a sequence $D_j $ of solutions of the isoperimetric problem with volumes 
  $Vol(D_j)\rightarrow 0$ for which $\partial D_j$ is not the graph on the sphere  
$\mathbb{S}^{n-1}$ of $T_p M$ of a function $u_j =r_j (1+v_j )$ where $||v_j ||_{C^{2,\alpha}}$ goes to $0$. This contradicts lemma \ref{tpreec}.
\end{Dem} 
\begin{Thm}\label{pblm1}
For all $n,r,Q,m,\alpha$ there exists $0<v_2=v_2(n,r,Q,m,\alpha)$ such that for all $M\in\mathcal{M}^{m,\alpha}(n,Q,r)$, $0<v<v_2$, if $D\subseteq M$ has volume $v$ and $I_M(v)=Area(\partial D)$ then $\partial D$ is a pseudo-bubble.
\end{Thm}
\begin{Dem}
   An analysis of the proof of theorem $1$ of \cite{NarAnn} shows how this application of the implicit function theorem gives a constant, say $C_0$ depending on $n,r,Q,m,\alpha$ such that the normal graph of a function $u$ on the unit tangent sphere centered at $p\in M$ with $||u||_{C^{2,\alpha}}\leq C_0$, solution of the pseudo-bubbles equation is of the form $\beta(p,r)$, $r<r_0$ then the argument given in theorem 3.1 of \cite{NarAnn} applies.           
\end{Dem}

\begin{Cor}\label{Kleiner1} Let $0<v<v_2$,  then for all $M\in\mathcal{M}^{m,\alpha}(n, Q, r)$, suppose that there exist a minimizing current $T$ for the isoperimetric problem with small enclosed volume $v$, $p\in M$ being its center of mass, and $St_p \leq Isom(M)$ being the stabilizer of $p$ for the canonical action of the group of isometries $Isom(M)$ of $M$. Then for all $k\in St_p$, we have $k(T)=T$.
\end{Cor}
\begin{Dem} Following theorem \ref{TT1}, $\partial T$ is the pseudo-bubble $\beta(p,r)$ where $\omega_n \rho^n =Vol(T)$. If $k\in St_p$, then, $k(\beta(p,r))=\beta(k(p),r*)$ for some small $r*$. For small volumes parameter $r$ is in one to one correspondence with parameter $v$, but $v$ is the enclosed volume and this does not change under the action of an isometry so by uniqueness of pseudo-bubbles we have that $r*=r$ hence $\beta(k(p),r)=\beta(p,r)$ and $k(T)=T$.  
\subsection{Proof of theorem \ref{1}.}
 For what follows it will be useful to give the definitions below. 
\begin{Def}
Let $(D_j)_j\subseteq\tau_M$ we say that $(D_j)_j$ is an \textbf{almost minimizing sequence in volume $v>0$} if \begin{enumerate}[(i):]
       \item $Vol(D_j)\rightarrow v$,
       \item $Area(\partial D_j)\rightarrow I_M(v)$.
\end{enumerate}
\end{Def}
\begin{Def}
Given $\phi:M\rightarrow N$ be a diffeomorphism between two Riemannian manifolds and $\varepsilon>0$. We say that $\phi$ is a \textbf{$(1+\varepsilon)$-isometry} if for every $x,y\in M$ holds $\frac{1}{1+\varepsilon}d_M(x,y)\leq d_N(\phi(x),\phi(y))\leq (1+\varepsilon)d_M(x,y)$.       
\end{Def}
For the convenience of the reader we have divided the proof into a sequence of lemmas. 
\end{Dem}
\subsubsection{Existence of a minimizer in a $C^{m,\alpha}$ limit manifold}
\begin{Lemme}\label{Isopcomparisoninfinity}
Let $M$ be with bounded sectional curvature and positive injectivity radius. $(M, p_j)\rightarrow (M_{\infty},p_{\infty})$ in $C^{m,\alpha}$ topology, $m\geq 1$. Then \begin{equation}\label{Isopcomparisoninfinity1}I_{M_{\infty}}\geq I_M.\end{equation}
\end{Lemme}
\begin{Dem} Fix $0<v<Vol(M)$. Let $D_{\infty}\subseteq M_{\infty}$ an arbitrary domain of volume $v=Vol_{g_{\infty}}(D_{\infty})$. Put $r:=d_H(D_{\infty},p_{\infty})$, where $d_H$ denotes the Hausdorff distance. Consider the sequence $\varphi_j:B(p_{\infty},r+1)\rightarrow M$, of $(1+\varepsilon_j)$-isometry given by the convergence of pointed manifolds, for some sequence $\varepsilon_j\searrow 0$. Set $D_j:=\varphi_{j}(D_{\infty})$ and $v_j:=Vol(D_j)$ it is easy to see that 
\begin{enumerate}[(i):]
          \item $v_j\rightarrow v$, 
          \item $Area_g(\partial D_j)\rightarrow Area_{g_{\infty}}(\partial D_{\infty})$.
\end{enumerate} 
(i)-(ii) are true because $\varphi_j$ is a $1+\varepsilon_j$ isometry. 

After this very general preliminary construction that doesn't requires any bounded geometry assumptions on $M$, except for the existence of a limit manifold along a sequence, we proceed to the proof of (\ref{Isopcomparisoninfinity1}) by contradiction. Suppose that there exist a volume $0<v<Vol(M)$ satisfying 
\begin{equation}\label{Isopcomparisoninfinity2}
          I_{M_{\infty}}(v)< I_M(v). 
\end{equation}
 Then there is a domain $D_{\infty}\subseteq M_{\infty}$ such that 
$$
I_{M_{\infty}}(v)\leq A_{g_{\infty}}(\partial D_{\infty})< I_M(v).
$$
As above we can find domains $D_j\subset M$ satisfying (i)-(ii). Unfortunately the volumes $v_j$ in general are not exactly equal to $v$. So we have to readjust the domains $D_j$ to get $v_j=v$, for every $j$, preserving the property $A_g(\partial D_j)\rightarrow A_{g_{\infty}}(\partial D_{\infty})$ as $j\rightarrow +\infty$, to get the desired contradiction. This can be done using the following construction that will be used in many places in the sequel. 
Looking carefully to the proof of the deformation lemma of \cite{RitGalli} and the compensation lemma of \cite{NarGD}, one can convince himself that it is possible construct domains $D^{\infty}_j\subseteq B(p_{\infty},r+1)\subseteq M_{\infty}$, as a small perturbation of $D_{\infty}$ such that   
\begin{equation}
    A_{g_{\infty}}(\partial D^{\infty}_j)\leq A_{g_{\infty}}(\partial D_{\infty})+c\tilde{v}_j,
\end{equation}
\begin{equation}
\tilde{v}_j\searrow 0,
\end{equation}
and 
\begin{equation}
Vol_g(\varphi_j(D^{\infty}_j))=v.
\end{equation}
The preceeding discussion show the existence of bounded finite perimeter sets (actually smooth domains) $D_j:=\varphi_j(D^{\infty}_j)\subset M$ satisfying the following properties
\begin{equation}
Vol_g(D_j)=v,
\end{equation}
\begin{equation}
 |A_g(\partial D_j)-A_{g_{\infty}}(\partial D^{\infty}_j)|\rightarrow 0,
\end{equation}
the last equation is true because $\varphi_j$ is a $1+\varepsilon_j$ isometry. In this way we obtain a sequence of domains $D_j$ having all volume $v$, such that 
\begin{equation}
A_g(\partial D_j)\rightarrow A_{g_{\infty}}(\partial D_{\infty})<I_M(v).
\end{equation}
The last equation is the desired contradiction, and the theorem follows from the arbitrariety of $v$. 
\end{Dem}

The next lemma is simply a restatement of theorem \ref{tprec1}.
\begin{Lemme}\label{restatementtprec1}
For all $n,r,Q,m,\alpha$, and $\varepsilon>0$ there exists\\ $0<v_3=v_3(n,r,Q,m,\alpha,\varepsilon)$ such that for all $M\in\mathcal{M}^{m,\alpha}(n, Q, r)$, there is a positive number $\eta=\eta(\varepsilon ,M)>0$ with the following properties\\ if $0<v=Vol(D)<v_3$, $\frac{Area(\partial D)}{I_M(Vol(D))}<1+\eta$ it follows that there exists $p=p_D\in M$, $R=C(n,r,Q,m,\alpha)v^{\frac{1}{n}}$ satisfying 
\begin{equation}
       \frac{Vol(D\Delta B(p,R))}{Vol(D)}\leq\varepsilon.
\end{equation}
\end{Lemme}
\begin{Dem}
As it is easy to check this lemma is a restatement of theorem \ref{tprec1} in an $\varepsilon$-$\delta$ language with a little extra effort about uniformity in the class $\mathcal{M}^{m,\alpha}(n, Q, r)$, after having observed that the constant $C$ used in the proof of lemma \ref{l2} depends only on $n,r,Q,m,\alpha$.     
\end{Dem}

\begin{Def} Let $M$ be a Riemannian manifold. $0<v<Vol(M)$ we say that $I_{M}(v)$ is \textbf{achieved} if there exists an integral current $D\subseteq M$ such that $Vol(D)=v$ and $Area(\partial D)=I_M(v)$.
\end{Def}
\begin{Lemme}\label{pbestimate}
For all $n,r,Q,m,\alpha$ there exist $0<v_4=v_4(n,r,Q,m,\alpha)$, $C_1=C_1(n,r,Q,m,\alpha)>0$ such that for all $M\in\mathcal{M}^{m,\alpha}(n,Q,r)$, $0<v<v_4$, with $I_M(v)$ achieved then 
\begin{equation} 
       I_M(v+h)\leq I_M(v)+C_1h v^{-\frac{1}{n}},
\end{equation}
provided that $v+h<v_4$.
\end{Lemme}
\begin{Dem} Let us define, $v_4=Min\{1, v_0, v_1, v_2\}$. Put 
$\psi_{M,p}(\tilde{v})=Area(\beta)^{\frac{n}{n-1}}$ where $\beta$ is the pseudo-bubble of $M$, centered at $p$ and enclosing volume $\tilde{v}$. Then $\tilde{v}\mapsto\psi_{M,p}(\tilde{v})$ is $C^1$ and $||\psi_{M,p}||_{C^1([0,v_4])}\leq C$ uniformly with respect to $M$ and $p$, i.e., $C=C(n,r,Q, m,\alpha)$, 
this is a nontrivial consequence of the proof of the existence of pseudo-bubbles that could be found in \cite{NarAnn}. When $v+h<v_4$, $$\psi_{M,p}(v+h)\leq\psi_{M,p}(v)+Ch.$$
\begin{equation}
      \begin{array}{lll}
    I_M(v+h) & \leq & \psi_{M,p}(v+h)^{\frac{n-1}{n}}\\
             & \leq & \psi_{M,p}(v)^{\frac{n-1}{n}}\left(1+\frac{Ch}{\psi_{M,p}(v)}\right)^{\frac{n-1}{n}}\\
             & \leq & \psi_{M,p}(v)^{\frac{n-1}{n}}\left(1+\frac{n-1}{n}C'h\right)\\
             & \leq & \psi_{M,p}(v)^{\frac{n-1}{n}}+C_1h v^{-\frac{1}{n}}\\
            & \leq & I_M(v)+C_1h v^{-\frac{1}{n}}.  
      \end{array}
\end{equation}
\end{Dem}


Now we want to apply the theory of convergence of manifolds suitably mixed with geometric measure theory to the isoperimetric problem for small volumes. Some parts of the proof are inspired from \cite{RRosales}, theorem 2.1.
\begin{Lemme}\label{elmpv}
For all $n,r,Q,m,\alpha$, there exists $0<v_6=v_6(n,r,Q,m,\alpha)$ such that for all $M\in\mathcal{M}^{m,\alpha}(n,Q,r)$, and for all $v$, with $0<v<v_6$ there is a sequence of points $p_j$, a limit manifold $(M_{\infty}, p_{\infty},g_{\infty})\in\mathcal{M}^{m,\alpha}(n,Q,r)$ and a domain $D_{\infty}\subset M_{\infty}$ such that 
\begin{enumerate}[(I):]
      \item \label{elmpvI} $(M,p_j,g)\rightarrow (M_{\infty}, p_{\infty},g_{\infty})$ in $C^{m,\beta}$ topology for $\beta<\alpha$,
      \item \label{elmpvIII}$I_{M_{\infty}}(v)=Area_{g_{\infty}}(\partial D_{\infty})$, so $I_{M_{\infty}}(v)$ is achieved, 
      \item \label{elmpvIV} $\partial D_{\infty}$ is a pseudo-bubble,
      \item \label{elmpvII} $I_M(v)=I_{M_{\infty}}(v)$.
\end{enumerate}            
\end{Lemme}
\begin{Dem} Fix $1>\delta>0$, and $\varepsilon>0$ such that 
            \begin{equation}\label{elmpv1bis}
                \frac{1}{2}\delta\frac{c_n}{C_1}>\gamma(\varepsilon)^{\frac{1}{n}}>0,
            \end{equation}
            with $\gamma=\gamma(\varepsilon)=\frac{\varepsilon}{1-\varepsilon}$. Observe that this is possible because $\gamma(\varepsilon)\rightarrow 0$ as $\varepsilon\rightarrow 0$.                
Set $v_6=Min\{v_0,v_2,v_3,v_4\}$ as obtained respectively in lemma \ref{BerMeyer}, \ref{restatementtprec1}, \ref{pbestimate} and theorem \ref{pblm1}. Let $0<v<v_6$. Let $D_j$ be a minimizing sequence in volume $v$ i.e. $Vol(D_j)=v$ and $Area(\partial D_j)\rightarrow I_M(v)$. Take now $j$ large enough to have $\frac{Area(\partial D_j)}{I_M(v)}<1+\eta_{\varepsilon}$ with $\eta_{\varepsilon}>0$ as in theorem \ref{restatementtprec1}. There exist $p_j$, $R$ s.t. 
$$\frac{Vol(D_j\Delta B(p_j,R))}{Vol(D_j)}\leq\varepsilon.$$
By theorem \ref{Fthmct} applied to the sequence of pointed manifolds $(M, p_j, g)_j\subset M^{m,\alpha}(n, Q, r)$ we obtain the existence of a pointed manifold $(M_{\infty}, p_{\infty}, g_{\infty})$ s.t. $(M, p_j)\rightarrow (M_{\infty}, p_{\infty}, g_{\infty})$ in $C^{4,\beta}$ topology. 

What we want to do in the sequel is to define domains $\tilde{D}_j^c\subseteq M_{\infty}$ (passing to a subsequence if necessary), that are images via the diffeomorphisms $F_j$ of $C^{4,\beta}$ convergence of a suitable truncation $D'_j$ of $D_j$ with balls whose radii $t_j$ are given by the coarea formula (because it is needed to control the amount of area added in the truncation procedure), to obtain an integral current $D_{\infty}\subseteq M_{\infty}$ s.t. $\tilde{D}_j^c\rightarrow D_{\infty}$ in $\mathcal{F}_{loc}(M_{\infty})$ topology. This goal will be achieved by taking an exhaustion of $M_{\infty}$ by geodesic balls, applying a standard compactness argument of geometric measure theory in each of these balls and using a diagonal process. 

Take a sequence of scales $(r_i)$, $i\geq 0$ satisfying $r_0\geq R$ and $r_{i+1}\geq r_i+2i$, consider an exhaustion of $M_{\infty}$ by balls of center $p_{\infty}$ and radius $r_i$, i.e. $M_{\infty}=\bigcup_i B(p_{\infty}, r_i)$. Then for every $i$ the convergence in $C^{4,\beta}$ topology gives existence of $\nu_{r_i}>0$ and diffeomorphisms $F_{j,r_i}:B(p_{\infty},r_i)\rightarrow B(p_j,r_i)$ for all $j\geq\nu_{r_i}$, that are $(1+\varepsilon_j)$-isometries for some sequence $0\leq\varepsilon_j\rightarrow 0$. 

At this stage we start the diagonal process, determining a suitable double sequence of cutting radii $t_{i,j}>0$ with $i\geq 1$ and $j\in S_i\subseteq\mathbb{N}$ for some sequence of infinite sets $S_1\supseteq ...\supseteq S_{i-1}\supseteq S_i\supseteq S_{i+1}\supseteq ...$, defined inductively. 
Before to proceed we recall the argument of coarea used in this proof repeatedly. For every domain $D\subseteq M$, every point $p\in M$, and interval $J\subseteq\mathbb{R}$ there exists $t\in J$ such that
\begin{equation}
  Area(D\cap(\partial B(p, t)))=\frac{1}{|J|}\int_J Area((\partial B(p, s))\cap D)ds\leq\frac{Vol(D)}{|J|}.
\end{equation}
We proceed as follow, cut by coarea with radii $t_{1,j}\in ]r_1, r_{1}+j[$ for $j\geq\nu_{r_2}$ we get domains $D'_{1,j}=D_j\cap B(p_j,t_{1,j})$, $D''_{1,j}=D_j-D'_{1,j}$ for $j$ large enough (i.e., $j\geq\nu_{r_1}$), satisfying 
\begin{equation}\label{elmpv1}
\left|Area(\partial D'_{1,j})+Area(\partial D''_{1,j})-Area(\partial D_j)\right|\leq\frac{v}{1}.
\end{equation}
It is worth to note here that (\ref{elmpv1}) is equivalent to 
$$2 Area(\partial D\cap B(p_j ; t_{1,j}))\leq v.$$
Consider the sequence of domains $\left(\tilde{D}_{1,j}=F_{j,r_2}^{-1}(D'_{1,j})\right)_j$ for $j\geq\nu_{r_2}$, it is true that 
\begin{enumerate}
      \item $Area(\partial D'_{1,j})\leq Area(\partial D_j)+2\frac{v}{1}\leq I_M(v)+2\frac{v}{1}$,
      \item $Vol(D'_{1,j})\leq v$, 
\end{enumerate}
so we have volume and boundary area, of the sequence of domains, bounded by a constant. A standard argument of geometric measure theory allows us to extract a subsequence $D'_{1,j}$ with $j\in S_1\subseteq\mathbb{N}$, converging on $B(p_{\infty},r_2)$ to a domain $D_{\infty,1}$ in $\mathcal{F}_{B(p_{\infty},r_2)}$. Now we look at the subsequence $D_j$ with $j\in S_1$ and repeat the preceding argument to obtain radii $t_{2,j}\in]r_2, r_3[$ and a subsequence $D'_{2,j}=D_{j}\cap B(p_j,t_{2,j})$ for $j\in S_1$ and $j\geq\nu_{r_3}$ such that 
\begin{equation}\label{elmpv2}
\left|Area(\partial D'_{2,j})+Area(\partial D''_{2,j})-Area(\partial D_{j})\right|\leq\frac{v}{2}.
\end{equation}
Analogously, the sequence $\left(\tilde{D}_{2,j}=F_{j,r_3}^{-1}(D'_{2,j})\right)_j$ for $j$ running in $S_1$ has bounded volume and bounded boundary area, so there is a convergent subsequence $\left(\tilde{D}_{2,j}\right)$ defined on some subset $S_2\subseteq S_1$ that is convergent on $B(p_{\infty},r_3)$ to a domain $D_{\infty,2}$ in $\mathcal{F}_{B(p_{\infty},r_3)}$. Continuing in this way, we obtain the existence of $S_1\supseteq ...\supseteq S_{i-1}\supseteq S_i$, radii
$t_{k,j}\in ]r_k, r_k+k[$, domains $D'_{i,j}=D_{j}\cap B(p_{j}, t_{i,j})$, $D''_{i,j}=D_{j}-D'_{i,j}$ satisfying 
\begin{equation}\label{elmpv3}
\left|Area(\partial D'_{kj})+Area(\partial D''_{kj})-Area(\partial D_j)\right|\leq\frac{v}{k},
\end{equation}
for all $1\leq k\leq i$ and $j\in S_k$ and for all $i\geq 1$. Moreover, putting $\tilde{D}_{k,j}=F_{j,r_{k+1}}^{-1}(D'_{k,j})$ for all $1\leq k\leq i$ and $j\in S_k$ we have convergence of $(\tilde{D}_{k,j})_{j\in S_k}$ on $B(p_{\infty},r_{k+1})$ to a domain $D_{\infty,k}$ in $\mathcal{F}_{B(p_{\infty},r_{k+1})}$ for all $i\geq 1$ and $k\leq i$. 
Let $j_i$ be chosen inductively so that 
\begin{equation}
j_i<j_{i+1}
\end{equation}
\begin{eqnarray}\label{elmpv4bis}       
       Vol(\tilde{D}_{i,\sigma_{i}(j_i)}\Delta D_{\infty,i})\leq\frac{1}{i},      
\end{eqnarray}
define $\sigma(i)=\sigma_{i}(j_i)$, then the sequence $\tilde{D}_i^c:=F_{\sigma(i),r_{i+1}}^{-1}(D'_{i,\sigma(i)})$ converges to $D_{\infty}=\bigcup_i D_{\infty,i}$ in $\mathcal{F}_{loc}(M_{\infty})$ topology.
Observe, here that $|t_{i+1}-t_i|>i$. From now on, we restrict our attention to the sequences $\bar{D}_i=D_{\sigma_i}$, $\bar{D}'_i=D'_{\sigma_i}$, $\bar{D}''_i=D''_{\sigma_i}$, then we will call always $D_i$, $D'_i$, and $D''_i$, by abuse of notation.  Put, also $F_i=F_{\sigma(i),r_{i+1}}$. Rename $i$ by $j$. From this construction we argue that passing possibly to a subsequence one can build a minimizing sequence $D_j$ with the following properties
\begin{enumerate}[(i):]
       \item \label{elmpv4} $\left|Area(\partial D'_j)+Area(\partial D''_j)-Area(\partial D_j)\right|\leq\frac{v}{j}$, for all $j$,
       \item \label{elmpv5}
$\lim_{j\rightarrow +\infty} Area_g(\partial D'_j)=\lim_{j\rightarrow +\infty} Area_{g_{\infty}}(\partial \tilde{D}^c_j)$,
       \item \label{elmpv6} $Vol(\tilde{D}^c_j)\rightarrow Vol(D_{\infty})=v_{\infty}$,  
       \item \label{elmpv7} $Area(\partial D_{\infty})\leq\liminf Area(\partial\tilde{D}^c_{j})$,             
       \item \label{elmpv8} $v\geq v_{\infty}\geq (1-\varepsilon)v>0$,       
       \item \label{elmpv9} $\frac{w_{\infty}}{v_{\infty}}\leq\gamma$ with $w_{\infty}=v-v_{\infty}$,
       \item \label{elmpv10} $I_{M_{\infty}}(v_{\infty})=Area(\partial D_{\infty})$,
       \item \label{elmpv10bis} $Area(\partial D_{\infty})=\liminf Area(\partial\tilde{D}^c_{j})$.
\end{enumerate}
(\ref{elmpv4}) follows directly by the construction of the sequences $(D'_j)$. (\ref{elmpv5}) is an easy consequences of the fact that the diffeomorphisms given by $C^{4,\beta}$ convergence are $(1+\varepsilon_j)$-isometry for some sequence $0\leq\varepsilon_j\rightarrow 0$. Let us denote 
$B_{r_{j+1}}=B(p_{\infty}; r_{j+1})$. To prove (\ref{elmpv6}) observe
\begin{eqnarray*}
|Vol(\tilde{D}_j^c)-Vol(D_{\infty})| & \leq & |Vol(\tilde{D}_j^c)-Vol(D_{\infty}\cap B_{r_{j+1}})| + Vol(D_{\infty}-B_{r_{j+1}})\\ 
                                     & \leq  & Vol((\tilde{D}_j^c\Delta D_{\infty})\cap B_{r_{j+1}}) + Vol(D_{\infty}-B_{r_{j+1}}),
\end{eqnarray*}
where $B_{r_{j+1}}$ is the ball $B(p_{\infty}, r_{j+1})$ in $M_{\infty}$. It follows now, by (\ref{elmpv4bis}) that $$\lim_{j\rightarrow\infty} Vol(\tilde{D}_j^c)=Vol(D_{\infty}).$$ On the other hand, the definition of the sets $\tilde{D}_j^c$ 
gives us $\{ D_j^c \}\rightarrow D$ in $\mathcal{F}_{loc}(M)$. Hence $Area(\partial D)\leq \liminf_{j\rightarrow\infty}Area(\partial\tilde{D}_j^c)$ by the lower semicontinuity of boundary area with respect to flat norm in $\mathcal{F}_{loc}(M)$ which actually proves (\ref{elmpv7}). In (\ref{elmpv8})
the first inequality is true because every $D_{\infty,i}$ is a limit in flat norm of a sequence of currents having volume less than $v$, the second beacuse the radii $r_i$ are greater than $R$ so $Vol(D_{\infty,i})\geq (1-\varepsilon)v$. (\ref{elmpv9}) follows easily by (\ref{elmpv8}). To show (\ref{elmpv10}) we proceed by contradiction. Suppose that there exists a domain $\tilde{E}\in\tau_{M_{\infty}}$ having $Vol(\tilde{E})=v_{\infty}$, $Area(\partial\tilde{E})<Area(\partial D_{\infty})$. Take the sequence of radii $s_j\in ]t_j, t_{j+1}[$ and cut $\tilde{E}$ by coarea obtaining $\tilde{E}_j:=\tilde{E}\cap B(p_{\infty}, s_j)$ in such a manner that  
\begin{equation}
       Area_{g_{\infty}}(\tilde{E}_j\cap\partial B(p_{\infty}, s_j))\leq\frac{v_{\infty}}{j},
\end{equation}
Of course, $Vol_{g_{\infty}}(\tilde{E}_j)\rightarrow v_{\infty}$, since $s_j\nearrow+\infty$.  
Now, fix a point $x_0\in\partial\tilde{E}$ and a small neighborhood $\mathcal{U}$ of $x_0$. For $j$ large enough $\mathcal{U}\subseteq B(p_{\infty}, r_j)$. Push forward $\tilde{E}_j$ in $M$ getting $E_{j}:=F_{j}(\tilde{E}_j)\subseteq B(p_j, r_{j+1})$ so readjusting volumes by modifying slightly $E_i$ in $F_i(\mathcal{U})$ contained in $B(p_j, t_{j+1})$, we obtain domains $E'_j\subseteq B(p_j, r_{j+1})$ with the properties 
\begin{equation}
      E'_j\cap D''_j=\emptyset, 
\end{equation}
\begin{equation}
      Vol_g(E'_j\cup D''_j)=v,
\end{equation}
\begin{equation}
      Area(\partial E'_j)\leq Area(\partial E_j)+c\Delta v_j,
\end{equation} 
with $\Delta v_j=Vol_g(E'_j)-Vol_g(E_j)$, satisfying $\Delta v_j\rightarrow 0$ as $j\rightarrow +\infty$, by virtue of $Vol(\tilde{E}_j)\rightarrow v_{\infty}$ (i.e. $Vol(D'_j)\rightarrow v_{\infty}$) and $Vol(D''_j)\rightarrow v-v_{\infty}$. Note that $c=c(n,Q)$ is a constant independent of $j$. Define $D^*_j:=E'_j\cup D''_j$. 
\begin{eqnarray*}
       Area(\partial D^*_j) & \leq & Area(\partial E'_j)+Area(D''_j)\\
            & \leq & (1+\varepsilon_j)^{n-1}Area(\partial\tilde{E}_j)+c\Delta v_j+Area(\partial D''_j)\\
    & \leq & (1+\varepsilon_j)^{n-1}(Area(\partial\tilde{E})+\frac{v_{\infty}}{j})+c\Delta v_j+Area(\partial D''_j),
\end{eqnarray*} 
hence we get
\begin{eqnarray*}
   \liminf_{j\rightarrow +\infty} Area(\partial D^*_j) & \leq & Area(\partial\tilde{E})+\liminf_{j\rightarrow +\infty} Area(\partial D''_j)\\
   & < & Area(D_{\infty})+\liminf_{j\rightarrow +\infty} Area(\partial D''_j )\\
   & \leq & \liminf_{j\rightarrow +\infty} Area(\partial D'_j)+\liminf_{j\rightarrow +\infty} Area(\partial D''_j)\\
   & \leq & I_M(v). 
\end{eqnarray*}
This means that the sequence of domains $D^*_j$ do better than the minimizing sequence $D_j$, which is a contradiction that proves (\ref{elmpv10}). The proof of (\ref{elmpv10bis}) is similar; in fact we only have to work with $D_{\infty}$ instead of $\tilde{E}$. We must remark that this can be done since the set of regular points in $\partial D_{\infty}\cap M_{\infty}$ is open, so inside it we can perform the suitable small deformation as in \cite{NarGD} Compesation Lemma and $\cite{RitGalli}$ Deformation Lemma.  Roughly speaking inside the regular part it is possible to make a smooth deformation of the domain at constant volume to produce a competitor with controlled area variation.  

Letting $i\rightarrow +\infty$ in (\ref{elmpv4}), taking into account (\ref{elmpv5}), (\ref{elmpv7}) and (\ref{elmpv10}), and Berard-Meyer inequality yields
\begin{equation}\label{elmpv11}
I_{M_{\infty}}(v_{\infty})+\delta c_nw_{\infty}^{\frac{n-1}{n}}\leq I_M(v).
\end{equation}  
It remains to prove that $v_{\infty}$ cannot be strictly less than $v$, by contradiction. We know that $v\leq v_4\leq v_2$ then $D_{\infty}$ is a pseudo-bubble as it is easy to check by corollary \ref{pblm1}. This allow one to have as a direct consequence of lemma \ref{pbestimate}, the following estimate
\begin{equation}\label{elmpv12}
   I_{M_{\infty}}(v)=I_{M_{\infty}}(v_{\infty}+w_{\infty})\leq I_{M_{\infty}}(v_{\infty})+C_1v_{\infty}^{-\frac{1}{n}}w_{\infty}.
\end{equation} 
Assume $w_{\infty}>0$. From (\ref{elmpv11}), (\ref{elmpv12}) and lemma \ref{Isopcomparisoninfinity} one deduce
\begin{equation}\label{elmpv13}
    I_{M_{\infty}}(v_{\infty})+\delta c_nw_{\infty}^{\frac{n-1}{n}}\leq I_M(v)\leq I_{M_{\infty}}(v)\leq I_{M_{\infty}}(v_{\infty})+C_1v_{\infty}^{-\frac{1}{n}}w_{\infty}.
\end{equation}
\begin{equation}\label{elmpv14}
    \delta c_nw_{\infty}^{\frac{n-1}{n}}\leq C_1v_{\infty}^{-\frac{1}{n}}w_{\infty}.
\end{equation} 
Dividing the above inequalities by $w_{\infty}^{\frac{n-1}{n}}$ and combining with (\ref{elmpv9}) we obtain  
\begin{equation}
      \gamma(\varepsilon)^{\frac{1}{n}}\geq\delta\frac{c_n}{C_1},   
\end{equation}
which by our choice of $\varepsilon>0$ contradicts (\ref{elmpv1bis}). So $w_{\infty}=0$, which means $v_{\infty}=v$ and clearly $I_{M_{\infty}}(v)=I_{M_{\infty}}(v_{\infty})$ which proves (\ref{elmpvIII}) and (\ref{elmpvIV}). To finish the proof, we need of a last argument that give us (\ref{elmpvII}). In fact 
\begin{eqnarray*}
 I_M(v) & = & \liminf Area(\partial D'_j)+\liminf Area(\partial D''_j)\\
   & = & I_{M_{\infty}}(v_{\infty})+\liminf Area(\partial D''_j)\\
   & = & I_{M_{\infty}}(v)+\liminf Area(\partial D''_j)\\
   & \geq & I_{M_{\infty}}(v).
\end{eqnarray*} 
Which combined with $I_M(v)\leq I_{M_{\infty}}(v)$ gives $I_M(v)=I_{M_{\infty}}(v)$ that is exactly (\ref{elmpvII}).\\
\textbf{Remark:} It is easy to check that $\liminf Area(\partial D''_j)=0$. 
\end{Dem} 
\paragraph{End of the proof of theorem \ref{1}.}

\begin{Dem}\\ Take $v*\leq v_6$. Suppose $0<v<v^*$. 

In first we show (\ref{MainI}) implies (\ref{MainII}).  Let $p_0$ be a point where $p\mapsto f(p,v)$ attains its minimum. We show by contradiction that $\beta(p_0,v)$ is a solution of the isoperimetric problem. Assume that there is no isoperimetric domain having volume $v$. Let $D_j$ be a minimizing sequence, $Vol(D_j)=v$, 
\begin{equation}\label{Main1} 
      Area(\partial D_j)\rightarrow I_M(v) <f_{M}(p_0,v)
\end{equation} and the isoperimetric profile is not achieved. The choice of $v^*$ ensures the existence of a pseudo-bubble $D_{\infty}\subseteq M_{\infty}$, and points $p_j$ satisfying (I)-(IV) of lemma \ref{elmpv}. Hence $I_M(v)=I_{M_{\infty}}(v)=Area(\partial D_{\infty})=f_{M_{\infty}}(p_{\infty}, v)$. A continuity argument with respect to $C^{4,\beta}$ convergence applies, giving $f_{M_{\infty}}(p_{\infty}, v)=\lim  f_{M}(p_j, v)$. Furthermore, since $p_0$ is a minimum point implies that $\forall j\:\: f_{M}(p_j,v)\geq f_{M} (p_0,v)$ from this one can argue finally that $f_{M_{\infty}}(p_{\infty}, v)\geq f_{M} (p_0,v)$ which contradicts (\ref{Main1}).\\ \indent
In second we show (\ref{MainII}) implies (\ref{MainI}). Let $D$ be an isoperimetric domain of sufficiently small volume, it follows from theorem
\ref{pblm1} that $D=\beta(p_0,v)$ for some point $p$ and small real $v$. This suffices to ensure that $p\mapsto f(p,v)$ attains its minimum at $p_0$. 

Finally, (\ref{Main0}) is a straightforward consequence of lemma \ref{elmpv}, noticing that for small volumes $I_M(v)=I_{M_{\infty}}(v)$ for some limit manifold $(M,\tilde{p}_{\infty},g_{\infty})$ obtained as the limit of the sequence $(M,\tilde{p}_j, g)$ for some sequence of points $\tilde{p}_j$. Furthermore, $I_{M_{\infty}}(v)=f_{M_{\infty}}(p_{\infty}, v)$ for some point $p_{\infty}$ possibly different from $\tilde{p}_{\infty}$. Now adjust the sequence of points $\tilde{p}_j$ to get a sequence of points $p_{j}\in M$ such that $(M,p_j, g)\rightarrow (M_{\infty},p_{\infty},g_{\infty})$ with the same $M_{\infty}$ as above. This goal could be achieved by taking as $p_j$ the points $p_j=F_j(p_{\infty})=F_{B_{M_{\infty}}(\tilde{p}_{\infty},R),j}(p_{\infty})$ for large $j$, where $R=d_{M_{\infty}}(\tilde{p}_{\infty},p_{\infty})+1$ and the $F_j$'s are the diffeomorphisms given by the $C^{m,\alpha}$ convergence.   
\end{Dem}
\section{Asymptotic expansion of the isoperimetric profile}\label{32}

We prove, now, theorem \ref{Cor1} stated in the introduction.

\begin{Dem} 
Let us just recall here the definition of $S=Sup_{p\in M}\{ Sc(p)\}$.
Let $(p_j)_j$ such that $Sc(p_j)\nearrow S$, take the sequence $(M,p_j,g)$ and apply theorem \ref{Fthmct} then we get the existence of $(M'_{\infty},p'_{\infty},g)$ such that passing to a subsequence, if needed, $(M,p_j,g)\rightarrow (M'_{\infty},p'_{\infty},g)$ in $C^{m,\beta}$ topology for $0<\beta<\alpha$. It is easy to check by a continuity argument that
\begin{equation}\label{Cor11}
       Sc_{M_{\infty}}(p'_{\infty})=S.
\end{equation}
From the definition of isoperimetric profile and lemma \ref{Isopcomparisoninfinity} follows 
\begin{equation}\label{Cor12}
      f_{M'_{\infty}}(p'_{\infty}, v)\geq I_{M_{\infty}}(v)\geq I_M(v).   
\end{equation}
Consider an arbitrary sequence of volumes $v_k\rightarrow 0$ and look at the corresponding $D_{v_k}$ we conclude that $$I_M(v_k)=I_{M_{\infty,k}}(v_k)=f_{M_{\infty, k}}(p_{\infty,k},v_k).$$
The sequence $(M_{\infty,k})$ belongs again to $\mathcal{M}^{4,\alpha}(n, Q, r)$ and an application of the fundamental theorem of convergence of manifolds to this sequence of manifolds produces a subsequence noted always with $v_k$, a limit manifold $(M_{\infty}, p_{\infty})$ with $(M_{\infty,k},p_{\infty,k})\rightarrow (M_{\infty},p_{\infty})$ in $C^{4,\beta}$ topology for every $0<\beta<\alpha$. From the latter construction it follows that
\begin{equation}\label{Cor13}
      I_M(v_k)\sim f_{M_{\infty}}(p_{\infty},v_k), k\rightarrow +\infty.
\end{equation} 
Combining (\ref{Cor11}), (\ref{Cor12}), (\ref{Cor13}), (\ref{Ipexp}) yields
\begin{equation}
    \frac{f_{M'_{\infty}}(p'_{\infty}, v_k)-c_nv_k^{\frac{n-1}{n}}}{v_k^{\frac{n+1}{n}}}\leq\frac{I_M(v_k)-c_nv_k^{\frac{n-1}{n}}}{v_k^{\frac{n+1}{n}}}   
\end{equation}     
From the asymptotic relation (\ref{Cor13}) letting $k\rightarrow +\infty$ we conclude that 
\begin{equation}
       -Sc_{M'_{\infty}}(p'_{\infty})\geq -Sc_{M_{\infty}}(p_{\infty}),
\end{equation}
that immediately gives
\begin{equation}\label{Cor14}
      S\leq Sc_{M_{\infty}}(p_{\infty}).
\end{equation}
Since the construction of $M_{\infty}$ permits us to have a sequence of points $p''_j\in M$ with $Sc_M(p''_j)\rightarrow Sc_{M_{\infty}}(p_{\infty})$ we obtain
\begin{equation}\label{Cor15}
      Sc_{M_{\infty}}(p_{\infty})\leq S. 
\end{equation}
(\ref{Cor14}), (\ref{Cor15}), and the arbitrarity of the sequence $v_k$, finally, give (\ref{fasymex}).
\end{Dem}
      \newpage
      \markboth{References}{References}
      \bibliographystyle{alpha}
      \bibliography{these}
      \addcontentsline{toc}{section}{\numberline{}References}
      \emph{Stefano Nardulli\\ Dipartimento di Metodi e Modelli Matematici\\ Viale delle Scienze Edificio 8 - 90128 Palermo\\ email: nardulli@unipa.it} 
\end{document}